\documentclass{ert-l}%
\issueinfo{00}
  {}
  {}
  {1995}
 
\copyrightinfo{1995}
  {American Mathematical Society}

\usepackage{amssymb}    
\usepackage{epsfig}     
\usepackage[all]{xy}

\newcommand{\CC}{{\mathbb C}}
\newcommand{\ZZ}{{\mathbb Z}}

\newcommand{\liea}[1]{\mathfrak{#1}}
\newcommand{\lieg}[1]{\mathrm{#1}}

\newcommand{\rk}{\operatorname{rk}}

\theoremstyle{plain}
\newtheorem{thm}{Theorem}[section]

\newtheorem{lemma}[thm]{Lemma}
\newtheorem{cor}[thm]{Corollary}

\newtheorem{prop}[thm]{Proposition}

\theoremstyle{definition}
\newtheorem{re}[thm]{Remark}

\newtheorem{defn}[thm]{Definition}
\newtheorem{recipe}[thm]{Recipe}

\begin{document}

\title{A normal form for admissible characters 
in the sense of Lynch}
\author{Karin Baur}
\thanks{Author supported by a DARPA grant and 
by Uarda Frutiger-Fonds (Freie Akademische Stiftung)}
\address{Karin Baur, Department of Mathematics, University of
California, San Diego, USA}
\email{kbaur@math.ucsd.edu}
\date{October 13, 2004}

\maketitle

\section*{Abstract}
Parabolic subalgebras $\liea{p}$ of 
semisimple Lie algebras 
define a $\ZZ$-grading of the Lie algebra. 
If there exists a nilpotent element in the 
first graded part of $\liea{g}$ 
on which the adjoint group of 
$\liea{p}$ acts with a dense orbit, the parabolic 
subalgebra is said to be nice. 
The corresponding nilpotent element is also 
called admissible. 
Nice parabolic subalgebras of simple Lie algebras 
have been classified. In the case of Borel subalgebras 
a Richardson element of $\liea{g}_1$ is exactly 
one that involves all simple root spaces. 
It is however difficult to write 
down such nilpotent elements for general 
parabolic subalgebras. In this paper we give 
an explicit construction of admissible 
elements in $\liea{g}_1$ that uses as few 
root spaces as possible. 
%

%
\section*{Introduction}
%
Let $\liea{g}$ be a semisimple Lie algebra over $\CC$, 
$\liea{p}\subset\liea{g}$ a parabolic subalgebra. 
There is a 
$\ZZ$-grade of $\liea{g}$, 
$\liea{g}=\sum_j\liea{g}_j$ such that 
$\liea{p}=\sum_{j\ge 0}\liea{g}_j$ and 
$\liea{n}:=\sum_{j>0}\liea{g}_j$ is 
the nilradical of $\liea{p}$. 
By a theorem of Richardson (\cite{ri}) there 
is always a Richardson element in $\liea{n}$, 
i.e. an element $X\in\liea{n}$ satisfying 
$[\liea{p},X]=\liea{n}$. 
We say that $\liea{p}$ is nice if there is a 
Richardson element in the first graded part 
$\liea{g}_1$. 
Nice parabolic subalgebras have been 
classified in~\cite{bw}. 
If $\liea{p}$ is a Borel subalgebra, then 
a Richardson element of $\liea{g}_1$ involves all 
simple root spaces.
For arbitrary nice parabolic subalgebras, the 
support of a Richardson element in the first 
graded part may consist of all roots of $\liea{g}_1$. 
In this sense it is far from being a simple 
representative of a Richardson element. 
The goal of this paper is to give a normal form 
of Richardson elements for nice parabolic subalgebras 
in the classical case. The construction 
uses as few root spaces of the nilradical as 
possible. It turns out that in many cases, the 
support of this normal form spans a simple system 
of roots. Since Richardson elements correspond 
to admissible characters 
$\nu:\overline{\liea{n}}\to\CC$, the normal form 
describes how admissible characters 
of the (opposite) nilradical actually look like.

In his thesis~\cite{l}, Lynch studied Whittaker 
modules for which there is an admissible homomorphism 
$\nu:\liea{n}\to\CC$: 
Let $\mathcal{U}$ be the universal enveloping algebra of 
$\liea{g}$ and let $V$ be a $\mathcal{U}$-module. 
A vector $v\in V$ is called a Whittaker vector if 
there exists a nonsingular character $\nu:\liea{n}\to\CC$ 
such that 
$xv=\nu(x)v$ for all $x\in\liea{n}$. The module $V$ 
is a Whittaker module if $V$ is cyclically generated 
by a Whittaker vector. 
An element $x\in\liea{g}$ is called admissible if 
$x\in\liea{g}_{-1}$ and $\liea{n}^x=\{0\}$. 
In particular, $x$ is a nilpotent element of $\liea{g}$. 
There is a natural bijection between $\liea{g}_{-1}$ 
and the characters of $\liea{n}$. That is 
if $x$ is an element of $\liea{g}_{-1}$, $x$ corresponds 
to the character $\nu:\liea{n}\to\CC$ if 
$B(x,y)=\nu(y)$ for all $y\in\liea{n}$. Here, 
$B(\cdot,\cdot)$ is the Killing form on $\liea{g}$. 
The character $\nu$ is called 
admissible if the corresponding element $x$ 
is admissible. 

Lynch has studied 
real parabolic subalgebras whose nilradicals 
support admissible Lie algebra homomorphisms to 
$i\mathbb{R}$. He calls a parabolic subalgebra 
$\liea{p}$ of $\liea{g}$ admissible if there exists 
is such an admissible homomorphism. 
Nice parabolic subalgebras are exactly the 
complexifications of admissible parabolic subalgebras 
in the sense of Lynch. 
So the classification of nice parabolic subalgebras 
tells us when admissible characters exist. 
But it is not easy to write one down. This paper 
gives a normal form for admissible elements of 
nice parabolic subalgebras in the classical 
Lie algebras. Via the correspondence 
admissible elements - admissible characters 
given above, the normal forms give a 
description of admissible characters. 

The basic idea of the construction is to spread 
out small blocks of 
skew-diagonal matrices at the appropriate places of 
the entries of the first graded part 
$\liea{g}_1$ of $\liea{p}$. 
In the fifth section we show that in many cases 
the normal form of a generic element 
involves only roots that form a simple system of roots. 
In particular, in all these cases we can construct 
a Richardson element using at most 
$n=\rk\liea{g}$ root spaces. 

The normal forms are given in section one. 
In section two we recall the properties 
of nice parabolic subalgebras. Sections three and four 
describe the constructions and give a modification 
for the orthogonal Lie algebras. 
In section five we discuss representation 
theoretic aspects of the results. We explain 
in which cases the roots involved in the normal form 
form a simple system of roots and describe the 
factors of the obtained root system. 
In the last section we turn our attention to 
the exceptional Lie algebras. Roughly half of the nice 
parabolic subalgebras of the exceptional Lie algebras 
have a Richardson element whose support forms a simple 
system of roots. We list all of these to complete the 
picture. 

Finally, I thank the referee for the suggestion concerning 
the exceptional Lie algebras. 

%
%
\section{Results}\label{section:results}
%
In what follows, $\liea{g}$ will be a classical 
Lie algebra over the complex numbers. 
As usual, 
$\liea{g}$ will be denoted by 
$\lieg{A_n},\lieg{B}_n,$ $\lieg{C}_n$, $\lieg{D}_n$ 
respectively. These Lie algebras are 
realized as subalgebras 
of the appropriate  $\liea{gl}_N$ ($N=n+1,2n+1,2n$ 
respectively $2n$). 
We follow the usage of~\cite{gw}: 
$\lieg{A}_n$ are the matrices of 
trace zero, $\lieg{B}_n$, $\lieg{D}_n$ the 
orthogonal Lie algebra of the symmetric form given 
by the matrix with entries $1$ on the 
skew-diagonal and otherwise $0$ (we will 
denote such a skew-diagonal $N\times N$-matrix by $J_N$). 
And 
$\lieg{C}_n$ is the symplectic Lie algebra 
for the symplectic form given by the matrix with 
entries $0$ outside the skew-diagonal and entry $1$ 
for the first $n$ entries on this skew-diagonal, 
entry $-1$ on the last $n$ entries of it (i.e. the form 
is given by the matrix $\begin{scriptsize} 
\begin{bmatrix}0 &J_n\\-J_n&0
\end{bmatrix}\end{scriptsize}$). 
Thus $\lieg{B}_n$ and $\lieg{D}_n$ are the 
matrices that are skew-symmetric around the skew-diagonal: 
If $A=(a_{ij})\in\liea{so}_N$ we have 
$a_{N-j+1,N-i+1}=-a_{ij}$. 
The elements of $\lieg{C}_n$ are the $2n\times 2n$-matrices 
$A=(a_{ij})$ with 
$a_{2n-j+1,2n-i+1}=-a_{ij}$ if $1\le i,j\le n$ and 
$a_{2n-j+1,2n-i+1}=a_{ij}$ if $1\le i\le n$, $n<j$ 
or $n<i$, $1\le j\le n$. 

With this realization the Borel subalgebra we choose is 
the intersection of the corresponding Lie algebra with 
the upper triangular matrices in $\liea{gl}_N$, the 
Cartan subalgebra is the set of diagonal matrices in 
the corresponding Lie algebra. We call a parabolic 
subalgebra containing this Borel subalgebra 
standard. 
Similarly, a Levi factor of $\liea{p}$ that contains 
the diagonal matrices is called a standard Levi 
factor of $\liea{p}$. 
If $\sum_j\liea{g}_j$ is the $\ZZ$-grade associated 
to $\liea{p}$ then $\liea{g}_0$ is the standard Levi 
factor of $\liea{p}$. 

The Levi factor $\liea{g}_0$ of a parabolic subalgebra 
consists of a sequence of square matrices on the diagonal. 
Caveat: in the case of $\liea{so}_N$, $\liea{sp}_{2n}$, 
the central block is skew-symmetric around the 
skew-diagonal or symmetric around the skew-diagonal, 
respectively. Thus in the case of $\lieg{D}_n$ 
there is an ambiguity in describing the parabolic 
subalgebras: 
the intersection of the standard parabolic 
subalgebra of $\liea{gl}_{2n}$ with block lengths 
$(a_1,\dots,a_{r-1},2,a_{r-1},\dots,a_1)$ with the 
$\liea{so}_{2n}$ is the same as the intersection 
of the parabolic subalgebra of $\liea{gl}_{2n}$ 
described by $(a_1,\dots,a_{r-1},1,1,a_{r-1},\dots,a_1)$ 
with $\liea{so}_{2n}$. In what follows we will 
always use the first version of this parabolic 
subalgebra of $\liea{so}_{2n}$. 

That said we can think of $\liea{p}$ to be given by the 
sequence of block lengths of the standard Levi factor. 
The space $\liea{g}_1$ is then a sequence of rectangles 
$R_{i,i+1}$ on the first super-diagonal. The sides of 
the rectangles are given by the lengths of the square 
matrices in $\liea{g}_0$. 
E.g. if $\liea{p}$ is given by $(a_1,\dots,a_r)$ then the 
$r-1$ rectangles of $\liea{g}_1$ 
have size 
$a_1\times a_2,a_2\times a_3,\dots,a_{r-1}\times a_r$. 

\begin{recipe}\label{recipe:A}
Let $\liea{g}$ be of type $\lieg{A}_n$, let the nice 
parabolic subalgebra 
$\liea{p}\subset\liea{g}$ be given by the unimodal sequence 
$(a_1,\dots,a_{r+1})$. Let 

$R_{i,i+1}$ have the form 
$\left\{\begin{array}{cl}
\begin{bmatrix} J_{a_i}& 0
\end{bmatrix} & \mbox{if $i$ is odd}, \\ 
 & \\
\begin{bmatrix}0 &J_{a_i} 
\end{bmatrix} & \mbox{if $i$ is even.}
\end{array}\right.$

Then we define $X_R\in\liea{g}_1$ to be the matrix 
formed by these rectangles. 
\end{recipe}
In the case of the symplectic and orthogonal Lie 
algebras with an even number of blocks in the standard 
Levi factor, $X_R$ is constructed similarly. 
However, if the standard Levi factor has an odd number 
of blocks we have to spread 
out the entries of the rectangles and choose small 
skew-diagonal matrices in both corners of 
each rectangle resp. in the central part of such 
a rectangle as we will see. To do so, set 
$B_i:=\lceil\frac{a_i}{2}\rceil$ and 
$b_i:=\lfloor\frac{a_i}{2}\rfloor$. 

If the standard Levi factor has an even number of 
blocks we talk about case (A), i.e. case (A) is the 
situation where the length in the standard Levi 
factor are given by $(a_1,\dots,a_r,a_r,\dots,a_1)$. 
If the lengths of the blocks in the standard Levi 
factor are $(a_1,\dots,a_r,a_{r+1},a_r,\dots,a_1)$ 
then we talk about case (B) (i.e. an odd number 
of blocks in the standard Levi factor).

\begin{recipe}\label{recipe:C}
Let $\liea{g}$ be of type $\lieg{C}$, 
$\liea{p}\subset\liea{g}$ a nice parabolic subalgebra, 
given by an unimodal sequence. 

Case (A). Let $R_r:=J_{a_r}$
and for $i\le r-1$ let 

$R_{i,i+1}$ have the form
$\left\{\begin{array}{cl}
\begin{bmatrix} J_{a_i}&0
\end{bmatrix} & \mbox{if $i$ is odd}, \\ 
 & \\
\begin{bmatrix} 0 &J_{a_i} 
\end{bmatrix} & \mbox{if $i$ is even.}
\end{array}\right.$

Case (B) 
For $i=1,\dots,r$ we let 

$R_{i,i+1}$ be of the form 
$\left\{\begin{array}{cl}
\begin{scriptsize}
\begin{bmatrix} &&J_{b_i} \\ J_{B_i}&&
\end{bmatrix}
\end{scriptsize} & \mbox{ if $i$ is odd,} \\ 
 & \\
\begin{scriptsize} 
\begin{bmatrix} &&J_{B_i} \\J_{b_i} &&
\end{bmatrix}\end{scriptsize} & 
\mbox {if $i$ is even.} 
\end{array}\right.$ 

Then we define $X_R\in\liea{g}_1$
to be the matrix that has the rectangles $R_{i,i+1}$ 
in the upper left half, that has the corresponding 
entries $\pm 1$ in the lower right half and - in case (A) 
of an even number of parts - 
has rectangle $R_r$ in the central position 
of the super-diagonal. 
\end{recipe}

In the cases of the orthogonal Lie algebras we allow 
the sequence of block lengths to be unimodal up to 
some smaller blocks in the middle. If $a_i>a_{i+1}$, 
the prescription is modified as follows. 
In case (A) 
we let $R_{i,i+1}$ have the form 
$\begin{scriptsize}
\begin{bmatrix}J_{a_{i+1}}\\ 0 \end{bmatrix}
\end{scriptsize}$ for odd $i$ 
resp. 
$\begin{scriptsize}
\begin{bmatrix} 0 \\ J_{a_{i+1}}\end{bmatrix}
\end{scriptsize}$ for even $i$.
In case (B) we set 
$R_{i,i+1}=
\begin{scriptsize}
\begin{bmatrix}  & J_{B_i}\\J_{B_i} \end{bmatrix}
\end{scriptsize}$. 
Note that the entries of the two $J_{B_i}$ overlap 
if $a_i$ is odd, that is the $B_i$th column has two 
nonzero entries. 

\begin{recipe}\label{recipe:BD}
Let $\liea{g}$ of type $\lieg{B,D}$, 
$\liea{p}\subset\liea{g}$ be a nice parabolic subalgebra, 
given by a unimodal sequence or let 
$a_1\le\dots\le a_l>a_{l+1}=\dots=a_r$ (resp. $=a_{r+1}$) 
with $a_{l+1}=a_l+1$. 

Case (A). 
For $i\le r-1$ with $a_i\le a_{i+1}$ 
let 

$R_{i,i+1}$ have the form 
$\left\{\begin{array}{cl}
\begin{bmatrix} J_{a_i}& 0
\end{bmatrix} & \mbox{if $i$ is odd}, \\ 
 & \\
\begin{bmatrix} 0&J_{a_i} 
\end{bmatrix} & \mbox{if $i$ is even.}
\end{array}\right.$

If $a_i>a_{i+1}$ let $R_{i,i+1}$ have the 
form 
$\begin{scriptsize}
\begin{bmatrix} 0 \\J_{a_{i+1}} \end{bmatrix}
\end{scriptsize}$ for odd $i$ 
resp. 
$\begin{scriptsize}
\begin{bmatrix} J_{a_{i+1}}\\ 0\end{bmatrix}
\end{scriptsize}$ for even $i$.

Furthermore we let $R_r$ be the matrix 
that has $b_r$ two-by-two 
blocks 
$\begin{scriptsize}\begin{bmatrix}
1&0\\0&-1
\end{bmatrix}\end{scriptsize}$ starting from the 
lower left corner if $r$ is odd resp. from 
the upper right corner if $r$ is even. 

Case (B). 
For $i\le r$ let $R_{i,i+1}$ have the form 

$\left\{ 
\begin{array}{cl}
\begin{bmatrix}  
0 & J_{B_i}& 0 
\end{bmatrix}
 & \mbox{if $a_i$, $a_{i+1}$ have the same parity,} \\ 
 & \\
\begin{scriptsize}
\begin{bmatrix} & J_{B_i-1}\\ & \\
J_{B_i}&
\end{bmatrix}
\end{scriptsize} & 
\mbox{if $a_i=a_{i+1}+1$ and $a_i$ is even, $i$ odd} \\
 & \\
\begin{scriptsize}
\begin{bmatrix} & J_{B_i}\\ & \\
J_{B_i-1}&
\end{bmatrix}
\end{scriptsize} & 
\mbox{if $a_i=a_{i+1}+1$ and $a_i$ is even, $i$ even} \\
 & \\
\begin{scriptsize}
\begin{bmatrix} & &J_{B_i}\\ 
J_{B_i}&&
\end{bmatrix}
\end{scriptsize} & \mbox{else (including $a_i=a_{i+1}+1$ 
with $a_i$ odd).}
\end{array}\right.$

Let $X_R\in\liea{g}_1$ be the matrix that has the 
rectangles $R_{i,i+1}$ in its upper left half 
and the corresponding entries 
$- 1$ in the lower right half and that has 
the square matrix $R_r$ as the $r$st rectangle 
in the super-diagonal in case (A). 
\end{recipe}

Note that in Recipe~\ref{recipe:BD} part (B), 
if $a_l=a_{l+1}+1$, with $a_l$ odd (i.e. 
the Lie algebra is an even orthogonal one) 
the recipe picks $a_l+1$ entries for the 
rectangle $R_{l,l+1}$. 
The rank of $R_{l,l+1}$ is $a_l-1$. 

\begin{thm}\label{thm:constr-ok}
Let $\liea{p}\subset\liea{g}$ be a nice 
parabolic subalgebra of a classical Lie algebra. 
Let $X_R$ be constructed as above 
(Recipes~\ref{recipe:A},~\ref{recipe:C},~\ref{recipe:BD}). 

Then $[\liea{p},X_R]=\liea{n}$. 
\end{thm}

\begin{re}
In their recent article~\cite{gr}, S. Goodwin and 
G. R\"ohrle use an alternative description 
of Richardson elements in $\liea{g}_1$ for 
parabolic subalgebras of $\liea{gl}_n$. They follow 
a construction of~\cite{bhrr}. 
The Richardson element obtained that 
way uses the identity matrix of size 
$\min\{a_i,a_{i+1}\}$ in the rectangle
$R_{i,i+1}$, starting in the upper left 
corner of the rectangle. This choice ensures 
that the matrix and all its powers have maximal rank. 

For our purposes the construction we present in 
this paper fits better since we are interested in 
the representation theoretic meaning of the root 
spaces involved. This will be discussed in 
Section~\ref{section:root-theory}.
\end{re}

%
\section{Background on nice parabolic subalgebras}
%

In this section we recall the properties of 
nice parabolic subalgebras that we will need 
later. The proofs of these statements 
can be found in~\cite{bw} and in~\cite{wa}.
The following result, Theorem 2.5 in~\cite{bw}, 
determines whether an 
element of $\liea{g}_1$ is a Richardson element or not. 

\begin{thm}\label{thm:dimension-levi}
Let $\liea{p}\subset\liea{g}$ be a parabolic subalgebra, 
$\liea{p}=\liea{m}\oplus\liea{u}$ where $\liea{m}$ 
is a Levi factor of $\liea{p}$ and $\liea{u}$ the 
corresponding unipotent radical. 
Let $\liea{g}=\sum\liea{g}_i$ be the corresponding 
grade. 
Then 
$x\in\liea{g}_1$ is a Richardson element of $\liea{p}$ 
if and only if $\dim\liea{g}^x=\dim\liea{m}$. 
\end{thm}

This dimension criterion is essential both in the 
classical and in the exceptional case. 
In the classical case, the method used in~\cite{bw} 
is the following: one calculates the Jordan normal 
form of a generic nilpotent element of $\liea{g}_1$. 
Then from the Jordan normal form, the dimension 
of a generic centralizer is computed. 
By Theorem~\ref{thm:dimension-levi}, the parabolic 
subalgebra is nice if and only this dimension is 
the same as the dimension of a Levi factor.

Next we recall the obtained characterization of nice 
parabolic subalgebras 
of the classical Lie algebras as it is given 
in the first section of~\cite{bw}: 
\begin{thm}\label{thm:typeA,C}
A parabolic subalgebra in a Lie algebra of a type 
$\lieg{A,C}$ 
is nice if and only if the corresponding sequence of 
block lengths in the standard Levi factor is unimodal. 
\end{thm}

\begin{thm}\label{thm:B}
A parabolic subalgebra in a Lie algebra of type 
$\lieg{B}_n$ is nice if and only if the corresponding 
sequence of block lengths in the standard Levi factor 
is unimodal or satisfies 
\[
a_1\le\dots< a_l>b_1=\dots=b_s
<a_l\ge\dots\ge a_1
\]
with $b_1=a_l-1$, $l\ge 1$.
\end{thm}

\begin{thm}\label{thm:D}
A parabolic subalgebra in a Lie algebra of type 
$\lieg{B}_n$ is nice if and only if the corresponding 
sequence of block lengths in the standard Levi factor 
has one of the following forms: 

1) It is unimodal with an odd number of blocks 

2) It is unimodal with an even number of blocks and odd 
block lengths occur exactly twice 

3) The block lengths satisfy 
\[
a_1\le\dots\le a_l>b_1=\dots=b_s<a_l\ge\dots\ge a_1
\]
with $b_1=a_l-1$, $a_l$ is odd and if $s$ is 
even, the odd block 
lengths occur exactly twice.
\end{thm}
%




To be able to compute the dimensions of the 
centralizer of a nilpotent element, we use the 
following well-known result (see~\cite{cm}): 

\begin{thm}\label{thm:dimension-centralizer}
Let $x\in\liea{g}$ be a nilpotent element, let 
$n_1\ge\dots\ge n_r>0$ be its Jordan normal 
form and $m_1\ge \dots\ge m_s>0$ the dual 
partition. 
Then the dimension of the centralizer of 
$x$ in $\liea{g}$ is 
\[\left\{ 
\begin{array}{ll}
\sum_i m_i^2 & \mbox{if}\ \liea{g}=\liea{sl}_{n+1} \\
\sum_i m_i^2+\frac{1}{2}\{i\mid n_i\ \mbox{odd}\}
 &\mbox{if}\  \liea{g}=\liea{sp}_{2n} \\
\sum_i m_i^2-\frac{1}{2}\{i\mid n_i\ \mbox{odd}\}
 &\mbox{if}\  \liea{g}=\liea{so}_n 
\end{array}\right.\]
\end{thm}

%
\section{Construction, Proofs}
%
%

Let $\liea{p}\subset\liea{p}$ be a nice parabolic subalgebra. 
Let the corresponding sequence of the block 
lengths of $\liea{g}_0$ be $(a_1,\dots,a_{r+1})$. 
Denote the entries of $X_R$ in the rectangles 
$R_{i,i+1}$ by $X_i$, $i=1,\dots,r$. 
In the other cases, we have (skew-)symmetries around 
the skew-diagonal. 
As before, case (A) is the case when the standard Levi 
factor has an odd number of blocks. In that case, 
$\liea{p}$ is given by 
$(a_1,\dots,a_r,a_r,\dots,a_1)$, the elements 
of $\liea{g}_1$ have block form, say
\[
X_R=\begin{scriptsize}
\left[\begin{array}
[c]{cccccccc}%
0 & X_{1} & 0 & 0 & 0 & 0 & 0 & 0\\
\vdots & \vdots & \ddots & \vdots & 
\vdots & \vdots & \vdots & \vdots\\
0 & 0 & 0 & X_{r-1} & 0 & 0 & 0 & 0\\
0 & 0 & 0 & 0 & X_r & 0 & 0 & 0\\
0 & 0 & 0 & 0 & 0 & Y_{r-1} & 0 & 0\\
\vdots & \vdots & \vdots & \vdots & 
\vdots & \vdots & \ddots & \vdots\\
0 & 0 & 0 & 0 & 0 & 0 & 0 & Y_{1}\\
0 & 0 & 0 & 0 & 0 & 0 & 0 & 0
\end{array}
\right]\end{scriptsize}
\]
If there is an even number of blocks, case (B), 
$\liea{p}$ is given by 
$(a_1,\dots,a_r,a_{r+1},a_r,\dots,a_1)$. The elements of 
$\liea{g}_1$ have block form, say 
\[
X_R=\begin{scriptsize}
\left[
\begin{array}
[c]{ccccccc}%
0 & X_{1} & 0 & 0 & 0 & 0 & 0\\
\vdots & \vdots & \ddots & \vdots & \vdots & \vdots & \vdots\\
0 & 0 & 0 & X_{r} & 0 & 0 & 0\\
0 & 0 & 0 & 0 & Y_{r} & 0 & 0\\
\vdots & \vdots & \vdots & \vdots & \vdots & \ddots & \vdots\\
0 & 0 & 0 & 0 & 0 & 0 & Y_{1}\\
0 & 0 & 0 & 0 & 0 & 0 & 0
\end{array}
\right]\end{scriptsize}
\]

Note that in the (A)-case $X_r$ is symmetric 
if the Lie algebra is 
symplectic, it is skew-symmetric if the Lie algebra is 
orthogonal. In the (A)-case, each $Y_i$ is  
skew-symmetric to $X_i$. 

In the (B)-case, if the Lie algebra is symplectic, then 
the first half of the rows of $Y_r$ are symmetric to the 
last half of the columns of $X_r$, and the last half of the 
rows of $Y_r$ are skew-symmetric to the first half of 
columns of $X_r$. For $i<r$, $Y_i$ is skew-symmetric to 
$X_i$. In the orthogonal case, each $X_i$ is skew-symmetric 
to $X_i$. 

\begin{re}\label{re:products-X_R}
Consider the powers $X_R^j$. 
The entries of $X_R^j$ are in the blocks in
the $j$-th diagonal, i.e. in the diagonal that consists 
of the rectangles 
$R_{i,i+j}$. 
\begin{enumerate}
\item
If $j=2$ with $a_i\le a_{i+1}\le a_{i+2}$, the products are 
described easily as of the form 
\[
\begin{bmatrix}
J_{a_i}&0\end{bmatrix}
\begin{bmatrix}
0&J_{a_{i+1}}\end{bmatrix}
=\begin{bmatrix}
0&\phantom{0}& I_{a_i}\end{bmatrix}
\]
or as 
\[
\begin{scriptsize}
\begin{bmatrix}
 & J_{B_i}\\
J_{b_i} & \end{bmatrix}
\begin{bmatrix}
 & J_{b_{i+1}} \\ 
J_{B_{i+1}}& \end{bmatrix}
=\begin{bmatrix}
I_{B_i} & \\ 
& I_{b_i} \end{bmatrix}
\end{scriptsize}
\]
or as 
\[
\begin{scriptsize}
\begin{bmatrix}
 & J_{B_i}\\
J_{B_i} & \end{bmatrix}
\begin{bmatrix}
 & J_{B_{i+1}} \\ 
J_{B_{i+1}}& \end{bmatrix}
=\begin{bmatrix}
I_{B_i} & \\ 
& I_{B_i} \end{bmatrix}
\end{scriptsize}
\]
depending on the case (where 
$b_i=\lfloor\frac{a_i}{2}\rfloor$, 
$B_i=\lceil\frac{a_i}{2}\rceil$).

\item More generally for $j\ge 1$ 
there are of three possible types:

\noindent1-A $X_{i}X_{i+1}\cdots X_{i+j-1}$ if $i\geq1$ and
$j+i\leq r+1$ (the $i$-th entry).

\noindent2-A $X_iX_{i+1}\cdots X_{r+1}Y_r\cdots Y_{2r-i-j+1}$
if $i\geq1$, $2r>i+j>r+1$ and $i\leq r$ (the $i$-th entry).

\noindent3-A $Y_{r-i}Y_{r-i-1}\cdots Y_{r-i-j+1}$ if
$i\geq0,i+j\leq r$ the ($r+i$-th entry).

In case (B) there are also three possibilities:

\noindent1-B $X_{i}X_{i+1}\cdots X_{i+j-1}$ if $i\geq1$ and
$j+i\leq r+1$ (the $i$-th entry).

%
\noindent2-B $X_{i}X_{i+1}\cdots X_{r}Y_{r}\cdots Y_{2r-i-j+1}$ if
$i\geq1$ and $2r>i+j>r+1$ and $i\leq r$ (the $i$-th entry).

\noindent3-B $Y_{r-i}Y_{r-i-1}\cdots Y_{r-i-j+1}$ if $i\geq0$ and
$i+j\leq r$ the ($r+i+1$-th entry).
\end{enumerate}
\end{re}

\begin{proof}[Proof of Theorem~\ref{thm:constr-ok}]
Let $\liea{p}\subset\liea{g}$ be a nice parabolic subalgebra, 
let $X_R\in\liea{g}_1$ be the corresponding nilpotent 
element. 

Since $\liea{p}$ is nice, there is a Richardson element 
$X$ in $\liea{g}_1$. By Theorem~\ref{thm:dimension-levi}, 
$\dim\liea{g}^X=\dim\liea{g}_0$. Recall that the 
dimension of the centralizer of a nilpotent element 
is calculated from its Jordan normal form (cf. 
Theorem~\ref{thm:dimension-centralizer}). The 
latter only depends on the ranks $\rk X^j$. 
So if we show that $X_R$ is generic (i.e. 
$\rk X^j=\rk X_R^j$ 
for all $j$), we obtain 
$\dim\liea{g}^X=\dim\liea{g}^{X_R}$. 
In particular, $X_R$ is a Richardson element. 

To see that $X_R$ is generic, we use 
the description of the blocks $X_R^j$ in 
Remark~\ref{re:products-X_R}. It is easy to 
check that in the case of $\lieg{A,C}$, 
the rank of $X_R^j$ is maximal for each $j$. 

In the $\lieg{B,D}$-case one has to be 
careful with the products of type 2-A and 2-B 
as described in Remark~\ref{re:products-X_R}. We 
leave the details to the reader and refer 
to~\cite{bw}, Section 3 where the generic 
ranks are computed explicitely in terms of the 
block lengths (Lemma 3.2 for $\lieg{A}_n$, 
Propositions 3.4 and 3.5 for the symplectic and 
the orthogonal Lie algebras).
\end{proof}

%
\section{More on the orthogonal Lie algebras}
%
We consider the case of a parabolic subalgebra 
$\liea{p}\subset\liea{g}$ of an orthogonal 
Lie algebra. The construction of the 
Richardson element $X_R$ given in 
Recipe~\ref{recipe:BD} seems unsatisfactory 
in the (B)-case, i.e. if the standard 
Levi factor has an odd number of blocks. 
In this section we will study this case in 
detail. The goal is to optimize the 
choice of $S_1$, the subset of the roots 
of $\liea{g}_1$ whenever this is possible. 
We will see that this is only possible if the sequence of 
the blocks in the standard Levi factor is not unimodal.

Let $\liea{p}$ be given by the sequence 
$(a_1,\dots,a_r,a_{r+1},a_r,\dots,a_1)$. Let 
$B_i:=\lceil a_i/2\rceil$ and $b_i=\lfloor a_i/2\rfloor$.  
Recall 
that the element $X_R$ has block form 
\[
X_R=\begin{scriptsize}
\left[
\begin{array}
[c]{ccccccc}%
0 & X_{1} & 0 & 0 & 0 & 0 & 0\\
\vdots & \vdots & \ddots & \vdots & \vdots & \vdots & \vdots\\
0 & 0 & 0 & X_{r} & 0 & 0 & 0\\
0 & 0 & 0 & 0 & Y_{r} & 0 & 0\\
\vdots & \vdots & \vdots & \vdots & \vdots & \ddots & \vdots\\
0 & 0 & 0 & 0 & 0 & 0 & Y_{1}\\
0 & 0 & 0 & 0 & 0 & 0 & 0
\end{array}
\right]\end{scriptsize}
\]
with 
$X_i=\left\{\begin{array}{cl}
(0\ \ J_{B_i}\ \ 0) & \mbox{if $a_i$ and $a_{i+1}$ 
have the same parity},\\
\begin{bmatrix} & J_{B_i} \\ & \\ J_{b_i} &\end{bmatrix} 
\mbox{or}
\begin{bmatrix} & J_{b_i} \\ & \\ J_{B_i} &\end{bmatrix} 
 & \mbox{if $a_i>a_{i+1}$, $a_i$ even, $i$ odd resp. $i$ even,} \\
\begin{bmatrix} & & J_{B_i}\\ J_{B_i}& &\end{bmatrix} 
 & \mbox{else}
\end{array}\right.$ 

\noindent 
and where $Y_i$ the negative of the skew-transpose of $X_i$.

In particular, if the parity of the $a_i$ changes from odd to 
even, the recipe chooses two roots of $R_{i,i+1}$ 
that subtract from each other. 
The next result shows that the extra root chosen is 
actually necessary in many cases. 

\begin{lemma}
Let the sequence $(a_1,\dots,a_{r+1},\dots,a_1)$ 
be unimodal, such that the corresponding parabolic 
subalgebra $\liea{p}\subset\liea{so}_n$ is nice. 
Assume that there is an $i\le r$ 
such that $a_i$ is odd and $a_{i+1}$ is even. 
Let $\hat{X}_i$ be a $a_i\times a_{i+1}$-matrix 
of rank $a_i$ 
with exactly $a_i$ nonzero entries, let $\hat{Y}_i$ 
be the negative of its skew-transpose.
Define $\hat{X}_R$ to be the matrix obtained from 
$X_R$ by replacing $X_i$ by $\hat{X}_i$ 
and $Y_i$ by $\hat{Y}_i$. 

Then $\hat{X}_R$ is not a Richardson element 
for $\liea{p}$. 
\end{lemma}
\begin{proof}
We have to show that $\hat{X}_R$ is not generic, in 
the sense that the ranks of $\hat{X}_R^j$ are 
not all maximal. 
By the construction of $X_R$, for each $j>i$ 
the matrix 
$Z_j:=X_j\cdots X_r\cdot Y_r\cdots Y_j$ is a 
$a_j\times a_j$ diagonal matrix of rank $a_j$ that 
is symmetric around the skew-diagonal. 

We can assume that $a_i=a_{i+1}-1$, the general statement 
follows with the same arguments. 
By assumption, $\hat{X}_i$ has a column of zeroes, say 
the $k$th. So the $(a_i+2-k)$th row of $\hat{Y}_i$ is 
zero. The product 
$\hat{X}_iZ_{i+1}\hat{Y}_i$ is a $a_i\times a_i$-matrix. 
Its entries are the entries of the diagonal matrix 
$Z_{i+1}$ except for the $k$th diagonal element of $Z_{i+1}$ 
and the $k$ to the last entry of $Z_{i+1}$. 
Since $a_{i+1}$ is even, these two entries of $Z_{i+1}$ are 
at different positions. Hence the product 
$\hat{X}_iZ_{i+1}\hat{Y}_i$ only has 
$a_{i+1}-2=a_i-1$ nonzero entries and is not 
of maximal rank.
\end{proof}
The next result treats the non-unimodal cases. 
Here we define $\hat{X}_i$ to be the matrix 
$\begin{scriptsize}\begin{bmatrix}
 & J_{b_i} \\J_{B_i} 
\end{bmatrix}\end{scriptsize}$ if $i$ is odd 
and 
$\begin{scriptsize}\begin{bmatrix}
 & J_{B_i} \\ J_{b_i} 
\end{bmatrix}\end{scriptsize}$ if $i$ is even.

\begin{lemma}\label{lemma:non-unim-B}
Let the sequence $(a_1,\dots,a_{r+1},\dots,a_1)$ 
describe the nice parabolic subalgebra 
$\liea{p}\subset\liea{so}_{2n+1}$. Assume that there is 
an $l$ such that $a_l=a_{l+1}+1$. Let $a_{l-1}$ be odd. 

Then the matrix $\hat{X}_R$ obtained by replacing 
$X_{l-1}$ ($Y_{l-1}$) by 
$\hat{X}_{l-1}$ (by $\hat{Y}_{l-1}$) 
is a Richardson element for $\liea{p}$ if and only if 
$a_{l-1}=a_{l+1}$. 
\end{lemma}

\begin{proof} 
As before we study the rank of the powers $\hat{X}_R^j$. 
Note that the matrices $X_{l+1},\dots,\ X_r,\ Y_r,\dots,Y_{r+1}$ 
all are square matrices with ones on the skew-diagonal. 
In particular, products of these are either equal to  
$I_{a_r}$ or to $J_{a_r}$. 
Hence we can assume w.l.o.g. that 
$l=r-1$, i.e. that the sequence satisfies 
\[
a_1\le \dots\le a_r>a_{r+1}<a_r\ge\dots\ge a_1.\] 

Let $a_{r-1}=a_{r+1}$, say $r-1$ is odd. Then 
\[\hat{X}_{r-1}X_r=
\begin{bmatrix}
 & J_{b_{r-1}}\\ J_{B_{r-1}}
\end{bmatrix}
\begin{bmatrix}
 & J_{B_{r+1}}\\ J_{b_{r+1}}
\end{bmatrix}
=\begin{bmatrix}
I_{b_{r+1}} & \\ & I_{B_{r+1}}
\end{bmatrix}
=I_{a_{r+1}}\]
which is equal to $-Y_r\hat{Y}_{r-1}$
and so 
$\hat{X}_{r-1}X_rY_r\hat{Y}_{r-1}=-I_{a_{r+1}}$. 

On the other hand, each product $X_k\cdots\hat{X}_i$ has 
(maximal) rank $a_k$ and the same holds for 
$\hat{Y}_i\cdots Y_k$. 
So in particular, all products of the blocks of $\hat{X}_R$ that 
involve $\hat{X}_i$ or $\hat{Y}_i$ have the same rank as 
the corresponding products for $X_R$. Therefore $\hat{X}_R$ 
is a Richardson element for $\liea{p}$. 

If $a_{r-1}\le a_{r_1}-2$, we obtain as 
\[\hat{X}_{r-1}X_r=
\begin{bmatrix}
 & J_{b_{r-1}}\\ J_{B_{r-1}}
\end{bmatrix}
\begin{bmatrix}
 & J_{B_{r+1}}\\ J_{b_{r+1}}
\end{bmatrix}
=\begin{bmatrix}
I_{b_{r-1}} && \\ && I_{B_{r-1}}
\end{bmatrix},
\] 
which is a matrix with $a_{r-1}$ rows and 
$a_{r+1}$ columns. So 
$Y_r\hat{Y}_{r-1}=-\begin{bmatrix}
I_{B_{r-1}}& \\ &  \\ & I_{b_{r-1}}
\end{bmatrix}$
and 
$\hat{X}_{r-1}X_rY_r\hat{Y}_{r-1}=
\begin{bmatrix}I_{b_{r-1}}& \\ & I_{b_{r-1}}
\end{bmatrix}$. 
Since $a_{r-1}$ is odd, the rank of this product 
is not maximal and therefore $\hat{X}_R$ is not 
a Richardson element of $\liea{p}$.
\end{proof}

We summarize these results and obtain one case where the Richardson 
element $X_R$ can be replaced by one with fewer root spaces involved. 

\begin{prop}\label{prop:BD-hat}
Let $\liea{p}\subset\liea{so}_N$ be a nice parabolic subalgebra 
with an odd number of blocks in the standard Levi factor. 
If $a_i$ is odd and $a_{i+1}$ is even then the element $\hat{X}_R$ obtained 
from $X_R$ by replacing $X_i$ and $Y_i$ by $\hat{X}_i$ and $\hat{Y}_i$ 
is a Richardson element for $\liea{p}$ exactly in the following case:

$\liea{g}=\liea{so}_{2n+1}$, $i=l-1$ with $a_{l-1}=a_{l+1}$ and 
$a_l>a_{l+1}$.
\end{prop}
From now on we will use the element $\hat{X}_R$ for non unimodal block lengths 
of $\liea{g}_0$ in $\liea{so}_{2n+1}$ whenever applicable. 
Let $S_1$ be the set of roots that are involved in the construction 
of the Richardson element. 
\begin{cor}\label{cor:non-simple}
Let $\liea{p}\subset\liea{so}_N$ be a nice parabolic 
subalgebra. Then there exist $\alpha,\beta$ in $S_1$ such 
that $\alpha-\beta$ 
is a root of $\liea{g}$ whenever 

(1) the sequence of the block lengths is unimodal, $a_i$ odd, 
$a_{i+1}$ even for some $i<r$.

(2) if $a_l>a_{l+1}$ 
there is $i<l$ with $a_i$ odd, $a_{i+1}$ even and if $i=l-1$ 
then $a_{r-1}<a_{r+1}$, 

(3) if $a_l>a_{l+1}$, $a_l$ is odd.
\end{cor}
\begin{proof}
In all these cases the recipes for $X_R$ resp. for $\hat{X}_R$ 
chooses more elements of $R_{i,i+1}$ than the ``rank'' of 
this rectangle, i.e. more than $\min\{a_i,a_{i+1}\}$.
\end{proof}

\begin{defn}\label{defn:star}
Let $\liea{p}\subset\liea{so}_N$ be a nice parabolic 
subalgebra as in (1), (2) or (3) of 
Corollary~\ref{cor:non-simple} above. Then we say 
that $\liea{p}$ is of the form (*).
\end{defn}
In other words, the parabolic subalgebras of the form 
(*) are exactly those where there are at least two roots in 
$S_1$ that subtract from each other. 
In the Section~\ref{section:root-theory} 
we will examine the structure of the roots of $S_1$ 
in detail.
%

%
%
\section{Root structure of $X_R$}\label{section:root-theory}
%

In this section we first translate the entries of $X_R$ to 
the corresponding roots. The constructed Richardson elements 
were chosen in a way that they involve as few root spaces as 
possible. We will see that very often the roots corresponding 
to $X_i$ do not subtract from each other. That means that 
they form a simple system of roots. We will discuss the different 
cases and show in which cases all roots of $X_R$ form a simple 
system. Also, we will explain what the factors of such a system 
are. 

Let $\liea{p}\subset\liea{g}$ be a parabolic 
subalgebra, let $\liea{m}=\liea{g}_0$ be its 
standard Levi factor. We denote the adjoint group 
of $\liea{m}$ by $M$. It acts on $\liea{g}_1$ 
by conjugation. Under this action, the rectangles 
$R_{i,i+1}$ (and there counterparts below the skew-diagonal 
in the case of $\lieg{B,C,D}$) 
are the irreducible components of 
$\liea{g}_1$. 

Note that $M$ is a product of the form 
$\lieg{GL}_{a_1}\times\cdots\times\lieg{GL}_{a_{r+1}}$ 
in the case of $\lieg{A}_n$. If 
$(h_1,\dots,h_{r+1})$ is an element of $M$, 
it sends $(X_1,\dots,X_r)\in\liea{g}_1$ 
to the element 
$(h_1^{-1}X_1h_2,h_2^{-1}X_2h_3,\dots,h_r^{-1}X_rh_{r+1})$. 

In the case of the other classical Lie algebras, 
(A)-cases, $M$ is of the form 
$\lieg{GL}_{a_1}\times\cdots\times\lieg{GL}_{a_r}$ 
and acts as 
$(h_1^{-1}X_1h_2,\dots,h_{r-1}^{-1}X_{r-1}h_r,h_r^{-1}X_r h_r)$.

In the (B)-cases, $M$ is the product 
$\lieg{GL}_{a_1}\times\cdots\times\lieg{H}_{a_r}$ 
with $H\subset\lieg{GL}_{a_{r+1}}$ the group 
$\lieg{SO}_{a_{r+1}}$ resp. $\lieg{SP}_{a_{r+1}}$. 
It acts on $\liea{g}_1$ as 
$(h_1^{-1}X_1h_2,\dots,h_{r-1}^{-1}X_{r-1}h_r,
h_r^{-1}X_r h_{r+1})$.

Under the action of the adjoint group $M$ 
the vector space $\liea{g}_1$ decomposes into 
$r$ irreducible components. 
The irreducible components are the vector spaces 
corresponding to the rectangles 
$R_{i,i+1}$ for $i\le r-1$ ($i\le r$ for $\lieg{A}_n$) 
together with $R_i$ 
in the (A)-cases (even number of blocks in the 
standard Levi factor). 
We thus obtain $r$ subsets of 
the positive roots of $\liea{g}$, in fact of 
the roots of $\liea{g}_1$. 
Each of these subsets fills out the corresponding 
rectangle $R_{i,i+1}$. 
The way we set up the construction of Richardson 
elements translates to the following: 

In the case of $\liea{sl}_{n+1}$ or in the (A)-cases 
(an even number of blocks in the standard Levi 
factor of 
the symplectic or of the orthogonal Lie algebra)
we choose $a_i$ entries starting 
alternatingly from the lower left corner 
resp. from the upper right corner. 
In terms of the roots of the rectangles, 
we choose the lowest resp. the highest root 
together with the next $a_i-1$ roots 
on the skew-diagonal. 
We describe the former case: 
In rectangle $R_{1,2}$ 
the roots are 
$\alpha_{a_1}$, $\alpha_{a_1-1,a_1,a_1+1}$, 
$\alpha_{a_1-2,\dots,a_1+2}$, up to 
$\alpha_{1,\dots,a_1,\dots,2a_1-1)}$ 
where $\alpha_{i,i+1,..,i+k}$ 
stands for the sum $\sum_{l=0}^k\alpha_{i+l}$.

In $R_r$ for $\liea{sp}_{2n}$ we choose all 
skew-diagonal entries. They correspond to 
the long roots of $\lieg{C}_n$ whose root 
spaces lie in $\liea{g}_1$. 

If the Lie algebra is $\liea{so}_{2n}$, 
the last rectangle $R_r$ is skew-symmetric 
around the skew-diagonal, so we start in  
the lower left or upper right corner just 
above the skew-diagonal: We choose 
the entries above (i.e. around) the skew-diagonal, 
leaving out every second. Say $r+1$ is odd. 
Then our choice corresponds to the roots 
$\alpha_{n-1}$, 
$\alpha_{n-3}+2\alpha_{n-2}+\alpha{n-1}+\alpha_n$, 
$\alpha_{n-5}+2\alpha_{n-4}+\dots+2\alpha_{n-2}
+\alpha_{n-1}+\alpha_n$, etc. 

(B)-cases. If there is an odd number of blocks in 
$\liea{g}_0$ we refined our choice of elements 
in the 
rectangles and started both from the upper 
right and the lower left corners. 
The pattern is similar as in case (A). 
The only difference is if for some $i$ 
$a_i$ is odd and $a_{i+1}$ is even, 
$\liea{g}=\liea{so}_N$. In that case we 
choose $a_i+1$ entries in $R_{i,i+1}$. 

To examine these roots better, let us denote 
the positive roots of $\liea{g}$ by $R^+$. 
Let the roots of the $i$th irreducible component 
of $\liea{g}_1$ under the adjoint group of the Levi 
factor be $R_1^i$. 
Let $S_1^i$ be the subsets of the roots whose 
root spaces contribute to $X_R$. 

On the matrix level, the roots in $S_1^i$ are  
the entries of $X_i$, $Y_i$ 
for $i<r$ and of $X_r$ in the (A)-case, 
they are the entries of $X_i$, $Y_i$ for 
$i\le r$ in the (B)-case 
and they are the entries of $X_i$ in the 
case of $\liea{sl}_{n+1}$. 

\begin{lemma}\label{lemma:root-difference}
Let $\alpha$, $\beta$ be elements of $S_1^i$, 
$\gamma\in S_1^j$. Then 

(1) $\alpha-\beta$ is not a  root of $\liea{g}$ except 
if $\liea{g}$ is an orthogonal Lie algebra and 
the parabolic subalgebra is of (*)-type. 

(2) 
$\alpha-\gamma$ is not a root of $\liea{g}$.
\end{lemma}

For the definition of the parabolic subalgebras of type 
(*) see Definition~\ref{defn:star}.
\begin{proof}
(1) 
Except for the (*)-cases, the recipes only 
choose roots that are in different rows/columns of 
$R_{i,i+1}$. 

(2) The roots of two different rectangles lie 
always in different rows. So in particular, their
difference is not a root of $\liea{g}$.
\end{proof}

\begin{lemma}\label{lemma:commute}
Let $\liea{p}\subset\liea{g}$ be a nice parabolic 
subalgebra with unimodal sequence of block 
lengths in $\liea{g}_0$. The roots of $S_1^i$ commute 
exactly in the following cases: 

(i) $\liea{g}$ is of type $\lieg{A}$ 
or $\lieg{C}$ case (A), $1\le i\le r$,

(ii) $\liea{g}$ is of type $\lieg{C}$ 
case (B) and $i<r$, 

(iii) $\liea{g}$ is of type $\liea{so}_{2n}$ case (A), 
$1\le i\le r$,

(iv) $\liea{g}=\liea{so}_N$ case (B), 
$i<r$ and if $a_i$ is odd there is no $k>i$ 
with $a_k$ even. 
\end{lemma}
\begin{proof}
In all these cases, the elements of $X_i$ (and of $Y_i$)  
are either diagonal neighbors (from lower left 
to upper right) or further apart from each other. 
If $\liea{g}$ is $\liea{sl}_{n+1}$ or if $i<r$, 
the roots of $X_i$ (and of $Y_i$) are the roots of 
$\liea{gl}_N$. Say we are dealing with entries (1,1) 
and (2,2) of the rectangle $R_{i,i+1}$. 
They correspond to the roots $\alpha_l$ 
and $\alpha_{l-1}+\alpha_l+\alpha_{l+1}$ for some 
$l$. Their sum is not a root of $\liea{g}$. 
In particular, 
$(\alpha_l,\alpha_{l-1}+\alpha_l+\alpha_{l+1})=0$. 

It remains to consider $i=r$, 
$\liea{g}=\lieg{C}_n$ or $=\lieg{D}_n$ 
with an even number of blocks in the standard 
Levi factor. In the former case, 
the recipe picks all long roots of $\liea{g}_1$. 
They do commute. 
In the latter case, the chosen entries of 
$R_{r,r+1}$ all contain $\alpha_n$. In particular, 
no two of them can add up to a root of $\liea{g}$. 
\end{proof}

The following result describes the structure of the 
roots of $S_1^r$ in the symplectic and orthogonal 
Lie algebras if the standard Levi factor of 
$\liea{p}$ has an odd number of blocks. 

\begin{lemma}
Let $\liea{p}$ be a nice 
parabolic subalgebra of $\liea{sp}_{2n}$ 
or of $\liea{so}_N$, given by the unimodal sequence 
$(a_1,\dots,a_r,a_{r+1},a_r,\dots,a_1)$. 

Then the roots in $S_1^r$ span 

$\left\{
\begin{array}{ll}
\frac{a_r}{2}\mbox{ factors }\lieg{A}_2 
& \mbox{if $a_r$ is even}, \\
\frac{a_r-1}{2}\mbox{ factors }\lieg{A}_2
\ \mbox{and one factor}\ 
\lieg{A}_1 & \mbox{if $a_r$ is odd 
and $\liea{g}=\liea{sp}_{2n}$}.
\end{array}\right.$
\end{lemma}

\begin{proof}
Recall that the roots chosen by the recipes form 
a rectangle of the form 
$\begin{scriptsize}\begin{bmatrix}
 & J_{B_r} \\ J_{b_r} & 
\end{bmatrix}\end{scriptsize}$ 
where $b_r=\lfloor a_r/2\rfloor$ and 
$B_r=\lceil a_r/2\rceil$. 
We consider the case of $\lieg{C}_n$, with 
sequence $(a_1,a_2,a_1)$. The case of the 
orthogonal Lie algebras is similar and thus 
is left to the reader. 

Let $1\le l\le a_1$. Then 
the $l$th entry from the lower left corner 
corresponds to the root 
$\alpha_{a_1-l+1}+\dots+\alpha_{a_1}+\dots+\alpha_{a_1+l-1}$. 
And the $l$th entry from the upper right corner 
corresponds to the root 
$\alpha_{l}+\dots+\alpha_{a_1}+\dots+\alpha_{a_1+l-1}
+2\alpha_{a_1+l}+\dots+2\alpha_{n-l}+\alpha_{n-l+1}$. 
Adding these two gives 
a root of $\liea{g}_2$. 
Then one checks that the $l$th entry from the lower left 
corner and the $k$th entry from the upper right 
corner do not add up to a root of $\liea{g}$ if 
$k\neq l$. 
\end{proof}

Now we have everything needed to understand the structure 
of $S_1:=\cup_i S_1^i$.
\begin{prop}
Let $\liea{p}\subset\liea{g}$ be a nice parabolic subalgebra, 
let $S_1$ be the roots that are involved 
in the construction of $X_R$.
The elements of $S_1$ form a simple system of roots 
exactly in the following cases: 

$\left\{\begin{array}{l}
\mbox{$\liea{g}$ is of type $\lieg{A}_n$ or $\lieg{C}_n$} \\  
\mbox{$\liea{g}$ is orthogonal and 
$\liea{p}$ is not as in (*)} 
\end{array}\right.$
\end{prop}

\begin{proof}
Lemma~\ref{lemma:root-difference} shows that no two roots 
of the set $S_1$ subtract from each other.
\end{proof}

\begin{cor}
Let $\liea{p}\subset\liea{g}$ be a nice parabolic 
subalgebra, let $\liea{g}$ have rank $n$. If 
$\liea{g}$ is a special or a symplectic Lie algebra, 
or if $\liea{p}\subset\liea{g}$ is not as in (*) then 
$S_1$ consists of at most $n$ roots.
\end{cor}

%
\section{The exceptional Lie algebras}
%

If $\liea{g}$ is one of the exceptional Lie algebras, 
we may also look for representatives of 
a Richardson element that involve as few root spaces 
as possible. It turns out that in many cases, the 
support of such an element does not form a simple 
system of 
roots.
We list all parabolic subalgebras of the exceptional 
Lie algebras where there actually exists a Richaradson 
element whose support forms a simple system 
of roots. In each of these cases we give an explicit 
representative. 
Thus we complete the construction of Richardson elements 
with simple support for the simple Lie algebras. 

We sketch the way to find the parabolic subalgebras 
with such a Richardson element and how to obtain 
a representative with simple support: 
By Theorem~\ref{thm:dimension-levi} 
we know that the dimension of a Richardson orbit 
${\mathcal O}_X$ for $\liea{p}$ is equal to 
$\dim\liea{g}-\dim\liea{m}$. 
In~\cite{bw}, all nice parabolic subalgebras of 
the exceptional Lie algebras are listed. 
For each of these, we compute the dimension of 
the Levi factor and thus obtain the dimension of 
the Richardson orbit. 
Using the lists of nilpotent orbits 
of~\cite{cm}, one finds the Bala-Carter labels 
for nilpotent orbits of the given dimension. 
The Bala-Carter label essentially gives the 
group for which 
the Richardson element is a regular nilpotent 
element. Note that in some cases there is 
ambiguity (i.e. there exist two or 
three orbits of the same dimension). 
An analysis of the roots of the irreducible 
components of $\liea{g}_1$ then helps to 
determine the correct label. 
Once the label is obtained, we have to look for 
such a simple system of roots among the roots 
of $\liea{g}_1$. 

We use the Bourbaki ordering of the 
simple roots (so in the case of $\lieg{E}_n$, 
$\alpha_2$ is the root at the top, attached 
to the third simple root of the horizontal string). 
For simplicity, we denote the positive roots 
of the exceptional 
Lie algebras by writing $\alpha_I$ where the 
index $I$ lists the simple roots involved. 
E.g. the root 
$\alpha_1+\alpha_2+\alpha_3+2\alpha_4+\alpha_5$ 
of $\lieg{E}_6$ will be written as 
$\alpha_{1234^25}$. The parabolic subalgebra will be 
described by an $n$-tuple $(u_1,\dots,u_n)\in\{0,1\}^n$: 
an entry $0$ at position $i$ means that the $i$th simple 
root is a root of the (standard) Levi factor, an entry $1$ at 
position $j$ means that the $j$th simple root is a root of 
the nilradical of $\liea{p}$ (as in~\cite{bw}). 

Note that there are two trivial cases: 
if all simple roots are roots of the Levi factor of $\liea{p}$ 
(i.e. $\liea{p}$ is given by 
the tuple $(0,\dots,0)$) then the Richardson 
orbit is the zero orbit. 
If the Levi factor of $\liea{p}$ is the Cartan 
subalgebra (i.e. $\liea{p}$ is given by the 
tuple $(1,\dots,1)$) then the nilpotent element 
$X=\sum_{\alpha\in\Delta} X_{\alpha}$ 
is a Richardson element (where $\Delta$ is 
the set of simple roots of $\liea{g}$).  
In $\lieg{G}_2$ the only parabolic  
subalgebras where the (minimal) support 
of a Richardson element forms a simple 
system of roots are the two trivial cases. 
So the first interesting case is $\lieg{F}_4$. 
In type $\lieg{F}_4$, we have several parabolic 
subalgebra where there exists a Richardson 
element with a simple support. 

In the following tables we list all nice parabolic subalgebras 
of the simple exceptional Lie algebras that have a Richardson 
elements whose support is a simple 
system of roots. The ones with a label $*$ to the left 
are corrected (there had been mistakes in the earlier 
version of this paper). 
For each of these we give the support 
of such a Richardson element 
$X_0$ in the second column. 

\[ 
\begin{array}{ccl}
\liea{p}\subset\liea{f}_4  & \mbox{Support of}\ X_0\\ 
\\
* & (1,0,0,0) &  
{\small
\xymatrix@-4mm{\alpha_{12^23^2}\ar@{-}[r] & \alpha_{123^24^2}
}} \\ \\
 & (0,0,0,1) &
{\small
\xymatrix@-5mm{
\alpha_{1234}\ar@{-}[r] & \alpha_{23^24}
}} \\ \\
 & (1,1,0,0) & 
{\small
\xymatrix@-4mm{
\alpha_{234}\ar@{-}[r] & \alpha_1\ar@2{->}[r] & \alpha_{23^2}
}}
\end{array}
\]
For $\lieg{E}_6$ we list the parabolic subalgebras up to the 
symmetry of the Dynkin diagram given by interchanging 
the vertices $1\leftrightarrow6$ and 
$3\leftrightarrow5$. 
\[ 
\begin{array}{ccll} 
& \liea{p}\subset\liea{e}_6&\phantom{--}  & \mbox{Support of $X_0$}\\ \\
 & (1, 1, 1, 0, 1, 0) & &
{\small
\xymatrix@-5mm{ & & \alpha_{56}\ar@{-}[d]  \\ 
 \alpha_{45}\ar@{-}[r]& \alpha_{2} \ar@{-}[r] & 
\alpha_{34}\ar@{-}[r] &\alpha_{1} 
}} \\
\\ 
 & (1, 1, 0, 1, 0, 1) & &
{\small
\xymatrix@-5mm{ & & \alpha_{56}\ar@{-}[d]  \\ 
 \alpha_{345}\ar@{-}[r]& \alpha_{2} \ar@{-}[r] & 
\alpha_{4}\ar@{-}[r] & \alpha_{13} 
}} \\ \\
 & (1, 1, 0, 0, 0, 1) & &
{\small
\xymatrix@-5mm{\alpha_{234}\ar@{-}[r] & \alpha_{456}\ar@{-}[r]
 & \alpha_{13}\ar@{-}[r]&\alpha_{245}}} \\ \\
* & (1, 1, 0, 0, 0, 0) &  &
{\small
\xymatrix@-5mm{\alpha_{2345}\ar@{-}[r] & \alpha_{134}\ar@{-}[r] 
& \alpha_{2456}}}\\ \\
 & (1, 0, 1, 1, 0, 1) & &
{\small
\xymatrix@-5mm{ & 
 & \alpha_{24}\ar@{-}[d] \\
\alpha_{6}\ar@{-}[r] & \alpha_{45}\ar@{-}[r] 
 & \alpha_{3}\ar@{-}[r]&\alpha_{1} }
}\\ \\
 & (1, 0, 1, 0, 0, 1) & &
{\small
\xymatrix@-5mm{\alpha_{345}\ar@{-}[r] & \alpha_{1}\ar@{-}[r] 
 & \alpha_{234}\ar@{-}[r] & \alpha_{56} & }}\\ \\
 & (1, 0, 1, 0, 0, 0) & &
{\small
\xymatrix@-5mm{\alpha_{234}\ar@{-}[r] 
&\alpha_{1}\ar@{-}[r]&\alpha_{345} & & &}
}\\ \\
 & (1, 0, 0, 1, 0, 0) & &
{\small
\xymatrix@-5mm{\alpha_{345}\ar@{-}[r] & \alpha_{1}\ar@{-}[r] 
 & \alpha_{234}\ar@{-}[r]&\alpha_{456} & \alpha_{245} }}
\\ \\
 & (1, 0, 0, 0, 0, 1) & &
{\small
\xymatrix@-5mm{\alpha_{1234}\ar@{-}[r]&\alpha_{3456} 
 &  \alpha_{1345}\ar@{-}[r]&\alpha_{2456}}
}
\\ \\
 & (1, 0, 0, 0, 0, 0) & &
{\small\xymatrix@-6mm{
\alpha_{1234} &  \alpha_{1345}
}}\\ \\
 & (0, 1, 0, 1, 0, 0) & &
{\small
\xymatrix@-5mm{ &  \alpha_{345}\ar@{-}[d] \\ 
  \alpha_{134}\ar@{-}[r] & \alpha_2\ar@{-}[r]&\alpha_{456} 
}} \\ \\
 & (0, 1, 0, 0, 0, 0) & &
{\small
\xymatrix@-5mm{
\alpha_{123^24^35^26}\ar@{-}[r] & \alpha_2
}}\\ \\
 & (0, 0, 1, 0, 0, 0) & &
{\small
\xymatrix@-5mm{
\alpha_{1234}\ar@{-}[r] & \alpha_{3456} & \alpha_{1345} 
 & \alpha_{2345}
}}
\end{array}\] 
\[ 
\begin{array}{ccll} 
& \liea{p}\subset\liea{e}_7& \phantom{--} & \mbox{Support of $X_0$}  \\ 
 &  \\ 
&(1, 0, 1, 1, 0, 1, 0) & & 
{\small
\xymatrix@-5mm{
 & & \alpha_1\ar@{-}[d] \\
\alpha_{56}\ar@{-}[r] & \alpha_4\ar@{-}[r]
 & \alpha_3\ar@{-}[r] & \alpha_{245}\ar@{-}[r] & \alpha_{67}
}} \\ \\
&(1, 0, 1, 0, 0, 1, 0) & & 
{\small
\xymatrix@-5mm{
 & & \alpha_{456}\ar@{-}[d] \\
\alpha_{345}\ar@{-}[r] & \alpha_1\ar@{-}[r] & 
\alpha_{234}\ar@{-}[r] & \alpha_{567} 
}} \\ \\
&(1, 0, 1, 0, 0, 0, 0)& & 
{\small
\xymatrix@-5mm{
 & \alpha_{23456}\ar@{-}[d] \\
\alpha_{234^25}\ar@{-}[r] & \alpha_1\ar@{-}[r] & \alpha_{34567} 
}} \\ \\
 & (1, 0, 0, 0, 0, 1, 1 )& & 
{\small
\xymatrix@-5mm{ 
\alpha_{1234}\ar@{-}[r] & \alpha_{3456}\ar@{-}[r] 
& \alpha_7\ar@{-}[r] & \alpha_{2456}\ar@{-}[r] & \alpha_{1345}
}} \\ \\
 & (1, 0, 0, 0, 0, 1, 0)& & 
{\small
\xymatrix@-5mm{ 
\alpha_{1234}\ar@{-}[r] & \alpha_{34567}\ar@{-}[r] 
 & \alpha_{2456}\ar@{-}[r] & \alpha_{1345}
}} \\ \\
* & (1, 0, 0, 0, 0, 0, 1)& & 
{\small
\xymatrix@-5mm{  
\alpha_{123456}\ar@{-}[r] & \alpha_{34567}\ar@{-}[r] 
 & \alpha_{1234^25} & \alpha_{24567}
}} \\ \\
 & (1, 0, 0, 0, 0, 0, 0)& & 
{\small
\xymatrix@-5mm{  
\alpha_1\ar@{-}[r] & \alpha_{12^23^34^45^36^27}
}} \\ \\
* & (0, 1, 0, 0, 0, 0, 0)& & 
{\small
\xymatrix@-5mm{  
\alpha_{234^25^26}\ar@{-}[r] & \alpha_{1234567} & \alpha_{1234^256}
 & \alpha_{234^2567} & \alpha_{123^24^25}
}} \\ \\
 & (0, 0, 1, 0, 0, 0, 1)& & 
{\small
\xymatrix@-5mm{  
\alpha_{2345}\ar@{-}[r] & \alpha_{4567}\ar@{-}[r]
 & \alpha_{1234}\ar@{-}[r] & \alpha_{3456} & \alpha_{1345}
}} \\ \\
 & (0, 0, 0, 1, 0, 1, 0)& & 
{\small
\xymatrix@-5mm{  
\alpha_{345}\ar@{-}[r] & \alpha_{67}\ar@{-}[r]
 & \alpha_{245}\ar@{-}[r] & \alpha_{134}\ar@{-}[r]
 & \alpha_{56}\ar@{-}[r] & \alpha_{234}
}} \\ \\
 & (0, 0, 0, 1, 0, 0, 0)& & 
{\small
\xymatrix@-5mm{  
\alpha_{3456}\ar@{-}[r] & \alpha_{1234}\ar@{-}[r]
 & \alpha_{4567}\ar@{-}[r] & \alpha_{2345} & \alpha_{1345}\ar@{-}[r] 
 & \alpha_{2456}
}} \\ \\
 & (0, 0, 0, 0, 1, 0, 1)& & 
{\small
\xymatrix@-5mm{  
\alpha_{1345}\ar@{-}[r] & \alpha_{2456}\ar@{-}[r]
 & \alpha_7\ar@{-}[r] & \alpha_{3456} & \alpha_{2345}
}} \\ \\
 & (0, 0, 0, 0, 1, 0, 0)& & 
{\small
\xymatrix@-5mm{  
\alpha_{123^24^2567}\ar@{-}[r] & \alpha_{2456}\ar@{-}[r] 
 & \alpha_{1345} & \alpha_{2345}\ar@{-}[r] & \alpha_{4567}
 & \alpha_{3456}
}} \\ \\
 & (0, 0, 0, 0, 0, 1, 0)& & 
{\small
\xymatrix@-5mm{  
\alpha_{56}\ar@{-}[r] & \alpha_{123^24^35^267}
 & \alpha_{67}\ar@{-}[r] & \alpha_{12^23^24^35^26}
}} \\ \\
 & (0, 0, 0, 0, 0, 0, 1) & & 
{\small
\xymatrix@-5mm{  
\alpha_7 & \alpha_{234^25^26^27} & \alpha_{1^22^23^34^45^36^27}
}} 
\end{array}\] 
\[ 
\begin{array}{ccll} 
 & \liea{p}\subset\liea{e}_8 & \phantom{-} 
& \mbox{Support of $X_0$}  \\  \\
 & (1, 0, 0, 0, 0, 1, 1, 1) &  & 
{\small
\xymatrix@-5mm{
 & & \alpha_8\ar@{-}[d] \\
\alpha_{1234}\ar@{-}[r] & \alpha_{3456}\ar@{-}[r] & \alpha_7\ar@{-}[r] 
 & \alpha_{2456}\ar@{-}[r] & \alpha_{1345}
}}\\ \\
 & (1, 0, 0, 0, 0, 1, 0, 0) & & 
{\small
\xymatrix@-5mm{  
\alpha_{34567}\ar@{-}[r] & \alpha_{1234}\ar@{-}[r]
 & \alpha_{45678}\ar@{-}[r] & \alpha_{23456}\ar@{-}[r] 
 & \alpha_{1345}\ar@{-}[r] & \alpha_{24567}
}}\\ \\
 & (1, 0, 0, 0, 0, 0, 1, 1)& & 
{\small
\xymatrix@-5mm{
 & & \alpha_{1234^25}\ar@{-}[d] \\
\alpha_{24567}\ar@{-}[r] & \alpha_8\ar@{-}[r]
 & \alpha_{34567}\ar@{-}[r] & \alpha_{123456}\\
}}\\ \\
* & (1, 0, 0, 0, 0, 0, 0, 1) & & 
{\small
\xymatrix@-5mm{ 
\alpha_{234^25^26^27^28}\ar@{-}[r] & \alpha_1\ar@{-}[r]
 & \alpha_{12^23^34^45^36^27}\ar@{-}[r] & \alpha_{8}
}}\\ \\
 & (1, 0, 0, 0, 0, 0, 0, 0)& & 
{\small
\xymatrix@-5mm{   
\alpha_{1234}\ar@{-}[r] & \alpha_{12^23^24^45^46^37^28} 
 & \alpha_{1345}\ar@{-}[r] & \alpha_{12^23^34^45^36^37^28}
}} \\ \\
 & (0, 1, 0, 0, 0, 0, 0, 1)& & 
{\small
\xymatrix@-5mm{
  & \alpha_{2345}\ar@{-}[d] \\
\alpha_{1234}\ar@{-}[r] & \alpha_{45678}\ar@{-}[r] 
 & \alpha_{123^24^25^26^27}
 & \alpha_{2456}\ar@{-}[r] & \alpha_{123^24^35^267} \\
}} \\ \\ 
* & (0, 0, 1, 0, 0, 0, 1, 0)& & 
\text{there is no Richardson element with simple support.}
\\ \\ 
 & (0, 0, 0, 0, 1, 0, 0, 1) & & 
{\small
\xymatrix@-5mm{
 & & \alpha_{1345}\ar@{-}[d] \\
\alpha_{3456}\ar@{-}[r] & \alpha_{78}\ar@{-}[r]
 & \alpha_{2456}\ar@{-}[r] & \alpha_{123^24^2567}
 & \alpha_{2345}\ar@{-}[r] & \alpha_{4567} \\
}} \\ \\ 
 & (0, 0, 0, 0, 0, 1, 0, 0) & & 
{\small
\xymatrix@-5mm{   
\alpha_{24567}\ar@{-}[r] & \alpha_{13456}\ar@{-}[r]
 & \alpha_{12^23^24^35^2678}\ar@{-}[r] & \alpha_{34567} 
 & \alpha_{23456}\ar@{-}[r] & \alpha_{45678} 
}}\\ \\
 & (0, 0, 0, 0, 0, 0, 1, 1)& & 
{\small
\xymatrix@-5mm{
 & \alpha_{123^24^2567}\ar@{-}[d] \\
\alpha_{234^25^26^27}\ar@{-}[r] & \alpha_8\ar@{-}[r] 
 & \alpha_{1234^25^267} 
}} \\ \\
 & (0, 0, 0, 0, 0, 0, 0, 1) & & 
{\small
\xymatrix@-5mm{ 
\alpha_8\ar@{-}[r] & \alpha_{1^22^33^44^65^56^47^38}
}}
\end{array}\] 


\end{document}